\documentclass{amsart}

\usepackage{amsmath, amssymb, amsfonts}
\usepackage[all]{xy}
\usepackage{enumitem}   
\usepackage{float}

\usepackage{multimedia}
\usepackage{pgf,pgfplots,pgfarrows,pgfnodes,pgfautomata,pgfheaps,pgfshade}
\usepackage{xmpmulti}
\usepackage{amsmath, amssymb, amsfonts}
\usepackage{ifsym}
\usepackage{epstopdf} 
\usepackage{epsfig}
\usepackage[normalem]{ulem}
\usepackage{ifthen}

\usepackage{tikz}
\usepackage{graphicx}
\usetikzlibrary{intersections,calc,positioning}
\usetikzlibrary{shapes.geometric,fit}
\usetikzlibrary{decorations,arrows,decorations.pathmorphing,backgrounds,positioning,fit,petri,shadows}
\usetikzlibrary{3d}

\newtheorem{tm}{Theorem}
\newtheorem{lm}[tm]{Lemma}
\newtheorem{prop}[tm]{Proposition}
\newtheorem{kor}[tm]{Corollary}

\newenvironment{dokaz}
{\noindent\emph{Proof:}\ }
{\hfill $\blacksquare$}

\newcommand{\Z}
{{\mathbb Z}}

\newcommand{\C}
{{\mathbb C}}

\newcommand{\g}
{{\mathfrak g}}

\newcommand{\gk}
{\hat{{\mathfrak g}}}
\newcommand{\gt}
{\tilde{{\mathfrak g}}}

\newcommand{\h}
{{\mathfrak h}}
\newcommand{\n}
{{\mathfrak n}}
\newcommand{\gsl}
{{\mathfrak sl}}

\newcommand{\Gamt}
{\tilde{{\Gamma}}}

\begin{document}

\author{Goran Trup\v{c}evi\'{c}}
\title{Bases of standard modules for affine Lie algebras of type $C_\ell^{(1)}$}
\address{University of Zagreb, Faculty of Teacher Education, Zagreb, Croatia}
\curraddr{}
\email{goran.trupcevic@ufzg.hr}
\thanks{Partially supported by the Croatian Science Foundation grant 2634.}
\date{}
\dedicatory{}

\begin{abstract}
Feigin-Stoyanovsky's type subspaces for affine Lie algebras of type $C_\ell^{(1)}$ have monomial bases with a nice combinatorial description. We describe bases of whole standard modules in terms of semi-infinite monomials obtained as ``a limit of translations'' of bases for Feigin-Stoyanovsky's type subspaces.
\end{abstract}

\maketitle

\section{Introduction}

Principal subspaces were introduced by B. Feigin and A. Stoyanovsky in 1994. 
who have recovered  Rogers-Ramanujan and Gordon
identities by computing character formulas for these spaces. Furthermore, from bases of these subspaces they have constructed bases of whole standard modules
 consisting of semi-infinite monomials (see \cite{FS}).

Let us briefly recall this construction of basis.
Let $\gk=\hat{\gsl}_2(\C)$ be an affine Lie algebra with basis $e(n)$, $h(n)$, $f(n)$; $n\in\Z$, and a central element $c$, and let $V=L(\Lambda_0)$ be the basic module for $\gk$, 
with a highest weight vector $v_0=v_{\Lambda_0}$.
Consider a commutative subalgebra
$\bar{\n}=\langle e(n) \rangle $.
Define the {\em principal subspace} of $V$ by $W=W(\Lambda_0)=U(\bar{\n})v_0$.
For a monomial $e(n_t) \cdots e(n_1)$ say that it satisfies {\em difference condition} (DC, for short) if $n_{s+1}\leq n_s-2$, $1\leq s <t$. Furthermore,
a monomial $e(n_t) \cdots e(n_1)$ satisfies {\em initial condition} (IC, for short) if $n_s<0$, $1\leq s \leq t$.
Then the set
$$\{e(n_t)\cdots e(n_1) v_{0}\,|\, \textrm{satisfies}\ DC\ \textrm{and}\ IC\}$$
is a basis of $W$.

	Set $v_{i}=T_i v_0$,  $W_i=T_i W_0=U(\bar{\n})v_i$, where $T_n,n\in\Z$ are translations from the affine Weyl group of $\gk$.  
	Then
\begin{equation}\label{W_incl}
\dots \subset W_1\subset W_0\subset W_{-1}\subset W_{-2}\subset \dots, \qquad V=\bigcup W_n,
\end{equation}
and 
	$$v_0=e(1)v_{-1} =e(1)e(3)v_{-2}=e(1)e(3)e(5)v_{-3}=\dots$$
Basis of $W_i$ consists of monomial vectors satisfying difference condition and a shifted initial condition: $n_s<-2i$.
The sequence of inclusions \eqref{W_incl} is described by
	$$e(n_t)\cdots e(n_1) v_{-i}=e(n_t)\cdots e(n_1) e(2i+1) v_{-i-1}=e(n_t)\cdots e(n_1) e(2i+1)e(2i+3) v_{-i-2}.$$
	If we now take ``the limit to $\infty$'' we obtain 
	$$v_{-i}=x_\alpha(2i+1) x_\alpha(2i+3) x_\alpha(2i+5)\cdots v_{-\infty}.$$
We say that a semi-infinite monomial $e(n_1)e(n_2)\cdots$ {\em stabilizes} if from some point on,
degrees of successive factors are successive odd numbers.	
	Then the set
	$$\{e(n_1)e(n_2)\cdots  v_{-\infty}\,|\, \textrm{satisfies}\ DC\ \textrm{and stabilizes}\}$$
	is a basis of $V$.
	
Similar in spirit to the notion of principal subspaces, {\em Feigin-Stoyanovsky's type subspaces} were 
introduced by M. Primc who has found bases of these subspaces, and in a similar fashion constructed from them bases of
whole standard modules. In \cite{P1} this was done for $\gt$ of type $A_\ell^{(1)}$, for all standard modules and a particular 
choice of Feigin-Stoyanovsky's type subspaces, and was carried further in \cite{P2} for all classical Lie algebras and all 
possible choices of Feigin-Stoyanovsky's type subspaces, but only for basic modules, and
in \cite{P3} for $\gt$ of type $B_2^{(1)}$, for all standard modules. 
In \cite{T1} and \cite{T2} we have constructed bases for $\gt$ of type $A_\ell^{(1)}$, for all standard modules and all possible 
choices of Feigin-Stoyanovsky's type subspaces. 

For $\gt$ of type $C_\ell^{(1)}$, a nice combinatorial description
of bases of Feigin-Stoyanovsky's type subspaces was given in \cite{BPT}. In this paper, we use the description of bases from \cite{BPT} 
to obtain bases of  whole standard modules consisting of semi-infinite monomials.

\section{Affine Lie algebra  $C_\ell^{(1)}$} \label{S: Afine}

Let ${\mathfrak g}$ be a complex simple  Lie algebra of type $C_\ell$. Let $\h$ be a Cartan
subalgebra of ${\mathfrak g}$ and ${\mathfrak g}={\mathfrak h}+\sum {\mathfrak g}_\alpha$ a root decomposition of ${\mathfrak g}$. Let 
$$
R=\{\pm \epsilon_i\pm \epsilon_j\,|\, 1\leq i\leq j\leq \ell\}\backslash\{0\}
$$ be the corresponding root system realized in $\mathbb R\sp\ell$ with
the canonical basis $\epsilon_1,\dots,\epsilon_\ell$.
Fix simple roots
$$\alpha_1=\epsilon_1-\epsilon_2,\quad\dots,\quad\alpha_{\ell-1}=\epsilon_{\ell-1}-\epsilon_\ell,\quad\alpha_\ell=2\epsilon_\ell$$ 
and let   $ \g=\n_-+ \h + \n_+$ be the corresponding triangular decomposition. Let
$\theta=2\alpha_1+\dots+2\alpha_{\ell-1}+\alpha_\ell=2\epsilon_1$ be the maximal root and
$$
\omega_r= \epsilon_1+\dots+ \epsilon_r, \qquad r=1,\dots,\ell
$$ 
fundamental weights (cf. \cite{H}). 
Fix root vectors $x_\alpha\in\g_\alpha$ and denote by $\alpha\sp\vee\in\mathfrak h$
dual roots. 
We identify $\mathfrak h$ and $\mathfrak h\sp*$ 
via the Killing form $\langle\,,\,\rangle$ normalized in such a way that  $\langle\theta,\theta\rangle=2$.

Denote by $\tilde{\mathfrak g}$ the affine Lie algebra of type $C_\ell\sp{(1)}$ associated  to $\g$,
$$
\hat{\mathfrak g} = \mathfrak{g} \otimes
\mathbb{C}[t,t^{-1}] + \mathbb{C}c, \qquad
\tilde{\mathfrak g} = \hat{\mathfrak g} + \mathbb{C} d,
$$
with the canonical central element $c$ and the degree element $d$ (cf. \cite{K}). It has a triangular decomposition 
$$
\gt=\tilde{\n}_-+ \tilde{\h} + \tilde{\n}_+,
$$
where $
\tilde{\n}_-=\n_-+ \g\otimes t^{-1}\C [t^{-1}]$, $\tilde{\h}=\h  + \mathbb{C}c + \mathbb{C} d$, $\tilde{\n}_+=\n_++ \g\otimes t\C [t]$.
Denote by $\Lambda_0,\dots,\Lambda_\ell$ fundamental weights of $\gt$.
 
For $x\in{\mathfrak g}$ and $n\in\mathbb Z$ denote by $x(n)=x\otimes t^{n}$ and $x(z)=\sum_{n\in\mathbb Z}
x(n) z^{-n-1}$, where $z$ is a  formal variable.



\section{Feigin-Stoyanovsky's type subspaces} \label{S: FS}

\begin{figure}[ht] \caption{Difference conditions} \label{figDC2}
	\begin{center}
		\begin{tikzpicture} [scale=.17]
		
		\begin{scope} [xshift=-40cm]
\node at (-6,16) {(a)};

		\draw (-.5,-.5) -- +(16,0) -- +(0,16) -- cycle;
		\node at (-1,15) {$\scriptscriptstyle 1$};
		\node at (-1,14) {$\scriptscriptstyle 2$};
		\node at (-1,0) {$\scriptscriptstyle \ell$};
		\node at (0,-1) {$\scriptscriptstyle 1$};
		\node at (1,-1) {$\scriptscriptstyle 2$};
		\node at (15,-1) {$\scriptscriptstyle \ell$};

		\node[fill=red, circle, inner sep=1pt](a) at (2,12) {};
		\node[fill=red, circle, inner sep=1pt](b) at (4,8) {};
		\node[fill=red, circle, inner sep=1pt](c) at (6,7) {};
		\draw[dotted] (a) -- (b) -- (c);
		
		\draw[dashed] (6,-.5) node[below]{$\scriptscriptstyle i_t$} -- (6,9) -- (-.5,9) ;

		\node[fill=blue, circle, inner sep=1pt](a1) at (1,8) {};
		\node[fill=blue, circle, inner sep=1pt](b1) at (5,7) {};
		\node[fill=blue, circle, inner sep=1pt](c1) at (7,3) {};
		\node[fill=blue, circle, inner sep=1pt](d1) at (9,2) {};
		\draw[dotted] (a1) -- (b1) -- (c1) -- (d1);
		\end{scope}
		
		\begin{scope}		
\node at (-6,16) {(b)};

		\draw (-.5,-.5) -- +(16,0) -- +(0,16) -- cycle;
		\node at (-1,15) {$\scriptscriptstyle 1$};
		\node at (-1,14) {$\scriptscriptstyle 2$};
		\node at (-1,0) {$\scriptscriptstyle \ell$};
		\node at (0,-1) {$\scriptscriptstyle 1$};
		\node at (1,-1) {$\scriptscriptstyle 2$};
		\node at (15,-1) {$\scriptscriptstyle \ell$};
		
		\node[anchor=-135] at (7.5,7.5) {$\scriptscriptstyle (-n)$};
		
		\node[fill=red, circle, inner sep=1pt](a) at (2,12) {};
		\node[fill=red, circle, inner sep=1pt](b) at (4,8) {};
		\node[fill=red, circle, inner sep=1pt](c) at (6,7) {};
		\draw[dotted] (a) -- (b) -- (c);
		\end{scope}
		\draw[dashed] (6,-.5) -- (6,9);
		
		\begin{scope}		
		\draw (-.5,-.5) -- +(16,0) -- +(16,-16) -- cycle;
		\node at (16,-1) {$\scriptscriptstyle 1$};
		\node at (16,-2) {$\scriptscriptstyle 2$};
		\node at (16,-16) {$\scriptscriptstyle \ell$};
		
		\draw[dashed] (6,-.5) -- (6,-7);
		\node[anchor=45] at (7.5,-8.5) {$\scriptscriptstyle (-n-1)$};
		
		\node[fill=blue, circle, inner sep=1pt](a1) at (7,-2) {};
		\node[fill=blue, circle, inner sep=1pt](b1) at (8,-6) {};
		\node[fill=blue, circle, inner sep=1pt](c1) at (12,-8) {};
		\node[fill=blue, circle, inner sep=1pt](d1) at (13,-10) {};
		\draw[dotted] (a1) -- (b1) -- (c1) -- (d1);
		\end{scope}
		\end{tikzpicture}
	\end{center}
\end{figure}
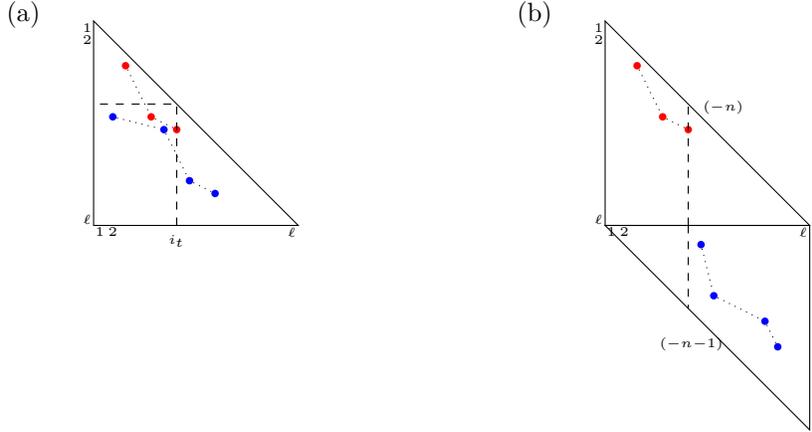

Fix the minuscule weight 
$\omega=\omega_\ell=\epsilon_1+\dots+\epsilon_\ell \in\h^*$; then $\langle\omega,\alpha\rangle\in\{-1,0,1\}$ for all $\alpha\in R$ and
define {\it the set of colors} 
$$\Gamma  =
\{\,\alpha \in R \mid \langle\omega,\alpha\rangle = 1\}  =  \{\epsilon_i+\epsilon_j\,|\, 1\leq i \leq j\leq \ell\}.
$$
Write
$$(ij)=\epsilon_i+\epsilon_j\in\Gamma\quad\textrm{and}\quad x_{ij}=x_{\epsilon_i+\epsilon_j}.$$
The set of colors can be pictured as a triangle with rows and columns ranging from $1$ to $\ell$ (see figure \ref{figDC2}, $(a)$); color $(ij)$ lies in the $i$-th column, and the $j$-th row of the triangle.

This gives a $\mathbb Z$-gradation of $\gt$; let ${\mathfrak g}_0 = {\mathfrak h} +
\sum_{\langle\omega,\alpha\rangle=0}\, {\mathfrak g}_\alpha$, then
$$
\gt=\gt_{- 1}+\gt_{0}+\gt_{ 1},
$$
where
$$
\gt_0 = \g_0 \otimes\C [t,t^{-1}] \oplus \C c \oplus \C d,\qquad
\gt_{\pm 1} = \sum_{\alpha \in \pm\Gamma}\, {\mathfrak g}_\alpha \otimes\C [t,t^{-1}].
$$
The subalgebra $\gt_1$ is commutative, and $\g_0$ acts on $\gt_1$ by adjoint action.

Let $L(\Lambda)$ be a standard $\tilde{\mathfrak g}$-module with the highest weight
$$
\Lambda=k_0 \Lambda_0+k_1 \Lambda_1+\dots+k_\ell \Lambda_\ell,
$$
$k_i\in\Z_+$ for $i=0,\dots,\ell$, and fix a highest weight vector $v_\Lambda$.
 Denote by
$k=\Lambda(c)$ the level of  $\tilde{\mathfrak g}$-module $L(\Lambda)$,
$
k=k_0 +k_1 +\dots+k_\ell$.

{\it Feigin-Stoyanovsky's type subspace} of $L(\Lambda)$ is  
$$
W(\Lambda)=U(\widetilde{{\mathfrak
		g}}_1)v_\Lambda\subset L(\Lambda).
$$

\bigskip

Before giving a description of bases of $W(\Lambda)$, we introduce a linear order on monomials from $U(\gt_1)$.
Define a linear order on $\Gamma$: set
$(i'j')<(ij)$ if $i'>i$ or $ i'=i,\, j'>j$. 
On the {\em set of variables} $\Gamt=\{x_\gamma(n) \,|\,\gamma\in\Gamma, n\in \Z \}$  set 
$x_\alpha(n)<x_\beta(n')$ if $n<n'$ or
$n=n',\, \alpha<\beta$. 
Unless otherwise specified, we assume that variables in monomials are sorted descendingly from right
to left. Define the order on the set of monomials as a lexicographic order -- compare variables from
right to left (from the greatest to the lowest one).

Order $<$ is compatible with multiplication (see \cite{P1}, \cite{T1}):
\begin{equation} \label{OrdMult_rel}
	\textrm{if} \quad x(\pi)<x(\pi')\quad \textrm{then}\quad x(\pi)x(\pi_1)<x(\pi')x(\pi_1)\quad 
\end{equation}
for monomials $x(\pi),x(\pi'),x(\pi_1)\in U(\gt_1)$.

\bigskip

\begin{figure}[ht] \caption{Basis of $W(\Lambda_r)$ consists of monomials satisfying difference and initial conditions} \label{figBasisFS}
	\begin{center}
		\begin{tikzpicture} [scale=.12]
		\clip (-3,16) rectangle (48,-29.5);
		
		\draw (-.5,-.5) -- +(16,0) -- +(0,16) -- +(0,0) -- +(16,-16) -- +(16,0) -- cycle;
		\draw (15.5,-16.5) -- +(16,0) -- +(0,16) -- +(0,0) -- +(16,-16) -- +(16,0) -- cycle;
		\draw (31.5,-32.5) -- +(16,0) -- +(0,16) -- cycle;
		\draw[dashed] (-.5,10) node[left] {$\scriptscriptstyle r$} -- (5,10);
		
		\node[anchor=-135] at (7.5,7.5) {$\scriptscriptstyle (-1)$};
		\node[anchor=45] at (7.5,-8.5) {$\scriptscriptstyle (-2)$};
		\node[anchor=-135] at (23.5,-8.5) {$\scriptscriptstyle (-3)$};
		\node[anchor=45] at (23.5,-24.5) {$\scriptscriptstyle (-4)$};
		\node[anchor=-135] at (39.5,-24.5) {$\scriptscriptstyle (-5)$};

		\node[fill=red, circle, inner sep=1pt](a) at (2,8) {};
		\node[fill=red, circle, inner sep=1pt](b) at (5,6) {};
		\node[fill=red, circle, inner sep=1pt](c) at (8,3) {};
		\node[fill=blue, circle, inner sep=1pt](d) at (10,-3) {};
		\node[fill=green, circle, inner sep=1pt](e) at (16,-5) {};
		\node[fill=green, circle, inner sep=1pt](f) at (18,-9) {};
		\node[fill=green, circle, inner sep=1pt](g) at (25,-10) {};
		\node[fill=orange, circle, inner sep=1pt](h) at (26,-18) {};
		\node[fill=orange, circle, inner sep=1pt](i) at (29,-23) {};
		\node[fill=yellow, circle, inner sep=1pt](j) at (34,-24) {};
		\node[fill=yellow, circle, inner sep=1pt](k) at (38,-25) {};
		\node[fill=yellow, circle, inner sep=1pt](l) at (42,-29) {};
		\draw[dotted] (a) -- (b) -- (c) -- (d) -- (e) -- (f) -- (g) -- (h) -- (i) -- (j) -- (k) -- (l);
		
		\end{tikzpicture}
	\end{center}
\end{figure}
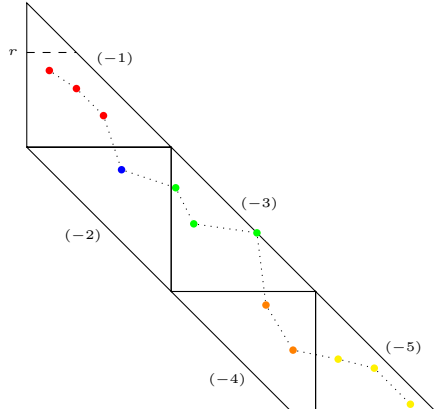

We say that a monomial $x(\pi)=\dots x_{i_s'
	j_s'}(-n-1)\cdots x_{i_1' j_1'}(-n-1)  x_{i_{t} j_{t}}(-n) \cdots x_{i_1
	j_1}(-n)\dots$ satisfies {\em
	difference conditions}, or shortly, that $x(\pi)$
satisfies {\em $DC$} for $W(\Lambda_r)$ if 
\begin{equation} \label{DClev1_eq}
i_1<\dots <i_t,\quad j_1<\dots<j_t,\quad  i_1'<\dots < i_s',\quad i_t<j_1'<\dots<j_s'. 
\end{equation}
This means that colors of factors of  $x(\pi)$ of each degree lie on diagonal paths in $\Gamma$, and a diagonal path of 
$(-n-1)$-part lies below $i_t$-th row, where $i_t$ is the column of the smallest
color of the $(-n)$-part; see figure \ref{figDC2}, (a).
Alternatively; below the triangle $\Gamma$ we can glue a transposed copy of $\Gamma$, with rows and columns interchanged, 
and represent the $(-n)$-part in the upper triangle and $(-n-1)$-part in the lower, transposed triangle. In this case, if 
a monomial satisfies difference conditions, then its $(-n-1)$-path lies on the right (and below) of the (-n) path, in the 
copy of $\Gamma^\tau$ (see figure \ref{figDC2}, (b))). If we continue this gluing and transposing procedure, we can represent
a monomial satisfying difference conditions by a diagonal path in a diagonal strip of copies of $\Gamma$ and $\Gamma^\tau$ 
(see figure \ref{figBasisFS}).

%
%
%

A monomial $x(\pi)$ satisfies {\em initial conditions} for $W(\Lambda_r)$ if 
it does not contain $x_{ij}(-1)$, $i\leq j \leq r$, as a factor. 
This means that the $(-1)$-path of $x(\pi)$ lies 
bellow the $r$-th row (see figure \ref{figBasisFS}).

Note that initial conditions can be interpreted in terms of difference conditions: 
$x(\pi)$ satisfies difference and initial conditions for $W(\Lambda_r)$ if and only if a monomial
$x(\pi')=x(\pi) x_{1r}(0)$ satisfies difference conditions.

In the level $k>1$ case, 
for $\Lambda=\Lambda_{r_1} + \dots +\Lambda_{r_k}$
we embed $L(\Lambda)$ in a tensor
product of standard modules of level $1$
$$L(\Lambda)\subset L({\Lambda_{r_1}}) \otimes \cdots \otimes L({\Lambda_{r_k}}),$$
with highest weight vector
$$v_\Lambda=v_{\Lambda_{r_1}} \otimes \cdots \otimes v_{\Lambda_{r_k}}.$$ 
In this case, a monomial
$x(\pi)$ satisfies {\em difference} and {\em initial conditions} for $W(\Lambda)$ if there exists a factorization
$x(\pi)=x(\pi^{(1)})\cdots x(\pi^{(k)})$ such that  $x(\pi^{(t)})$ satisfies difference and initial
conditions for $W(\Lambda_{r_t})$. 


\begin{tm}[{\cite{BPT}}]\label{FSbaza_tm}
	The set
	$$\{x(\pi)v_\Lambda \,|\, x(\pi) \ \textrm{satisfies}\ DC\ \textrm{and}\ IC\ \textrm{for}\ W(\Lambda)\}$$
	is a basis of $W(\Lambda)$.
\end{tm}

\section{Bases of standard modules} \label{S: Afine}

For a root $\alpha\in R$, let $x_\alpha,x_{-\alpha}$ be chosen such that $[x_\alpha,x_{-\alpha}]=-\alpha^\vee$. 
Define a ``Weyl group translation'' operator (cf. \cite{K}, \cite{FK})
\begin{eqnarray*}
	s_\alpha & = & \exp x_\alpha(0) \exp x_{-\alpha}(0) \exp x_\alpha(0),\\
	s_{\delta - \alpha} & = & \exp x_{-\alpha}(1) \exp x_{\alpha}(-1) \exp x_{-\alpha}(1),\\
	e_\alpha & = & s_{\delta - \alpha} s_\alpha.
\end{eqnarray*}
on $L(\Lambda)$. 
The following commutation relations hold
\begin{eqnarray}
\label{ekom1} e_\alpha d e_\alpha^{-1} & = & d + \alpha^\vee -\frac{1}{2}\langle \alpha^\vee,\alpha^\vee \rangle, \\
e_\alpha h e_\alpha^{-1} & = & h - \langle \alpha^\vee, h \rangle c, \\
e_\alpha h(n) e_\alpha^{-1} & = & h(n) \qquad \textrm{for}\ j\neq 0, \\
\label{ekom2} e_\alpha x_\beta(j) e_\alpha^{-1} & = & (-1)^{\beta(\alpha^\vee)} x_\beta(n-\beta(\alpha^\vee)),
\end{eqnarray}
for $h\in\h$, $\beta\in R$, $n\in\Z$.
Set
$$e=\prod_{\alpha\in\Gamma} e_\alpha.$$

The following proposition was proven by Primc in the case of $A_\ell$ (cf.
\cite{P1}, Theorem 8.2; \cite{P2}, Proposition 5.2)  and in the case of $B_2$(\cite{P3}, Proposition 8.2), 
and it carries over into this case without significant changes:

\begin{prop}[\cite{P1},\cite{P2},\cite{P3}] \label{StModGenSkup_prop} 
	$L(\Lambda)=\langle e \rangle U(\gt_1) v_\Lambda$.
\end{prop}

The main element of the proof is a generalization of Frenkel-Kac vertex operator formula (cf. \cite{FK}, \cite{F}, \cite{LP}, \cite{P1},\cite{P2},\cite{P3})
for $\hat{\gsl}_2(\alpha)$ generated by $x_\alpha(n)$, $x_{-\alpha}(n)$, $\alpha^\vee(n)$; $n\in\Z$, and a canonical central element $c_\alpha=2c/\langle \alpha,\alpha \rangle$, 
$$ \exp (zx_\alpha(z)) = E^-(-\alpha,z) \exp (-zx_{-\alpha}(z)) E^+(-\alpha,z) e_\alpha z^{c_\alpha+\alpha^\vee}.
$$
Note that for a standard module $L(\Lambda)$ of level $k$, the restriction to $\hat{\gsl}_2(\alpha)$ is of level $2k$ for $\alpha$ a short root, and of level $k$ for $\alpha$ a long root.

\bigskip

Basic module $L(\Lambda_0)$ is a vertex operator algebra, and $L(\Lambda_1),\dots,L(\Lambda_\ell)$ are modules for this algebra (cf. \cite{LL}). In that setting, operator $e$ can be described in terms of simple current operators. Let 
$$L(\Lambda_r) \xrightarrow{[\omega]} L(\Lambda_{\ell-r}) \xrightarrow{[\omega]}
L(\Lambda_i) $$
be simple current operators on level $1$ standard modules, such that the commutation relation 
\begin{equation} \label{OmegaComm} 
x_\alpha(n)[\omega]=[\omega] x_\alpha (n+\alpha(\omega)), \qquad \alpha\in R, \quad n\in\mathbb Z,
\end{equation}
holds (cf. \cite{DLM} and \cite{Ga}, also see Remark 5.1 in \cite{P3}). 

Then
$$	x_\gamma(-n-1) [\omega]  =  [\omega] x_\gamma(-n)$$
for $\gamma \in \Gamma$. Denote by $x(\pi^\pm)$ a monomial obtained by raising/decreasing degrees in $x(\pi)$ by $1$. By
$x(\pi^{\pm m})$ denote a monomial obtained by raising/decreasing degrees in $x(\pi)$ by $m$. Then
$$
	x(\mu)[\omega]  =  [\omega]x(\mu^+) 
$$
for a monomial $x(\mu)\in U(\gt_1)$.

In the level $k>1$ case, 
$L(\Lambda)\subset L(\Lambda_{0})^{k_0}\otimes
L(\Lambda_{1})^{k_1}\otimes \cdots \otimes
L(\Lambda_{\ell})^{k_\ell},$
use tensor products of level $1$ simple current operators
$$[\omega]=[\omega]\otimes\dots\otimes [\omega].$$

Like in \cite{P3},
\begin{lm}
	$e= C [\omega]^{2\ell}$, for some $C\neq 0$.
\end{lm}
\begin{dokaz}
	Note that
	\begin{eqnarray*}
		\sum_{\gamma\in\Gamma}\gamma^\vee & = & \sum_{1\leq i\leq j \leq \ell} (\epsilon_i+\epsilon_j)^\vee \\
		& = & \sum_{1\leq i < j \leq \ell} 2(\epsilon_i+\epsilon_j) + \sum_{i=1}^{\ell} 2\epsilon_i\\ 
		& = & 2 \ell (\epsilon_1+\dots+\epsilon_\ell)\\ 
		& = & 2 \ell \omega.
	\end{eqnarray*}
	Relations \eqref{ekom1}--\eqref{ekom2} and \eqref{OmegaComm} give 
	$e x_{\pm \gamma}(j) e^{-1} = x_{\pm \gamma}(j \pm 2 \ell)$ and $[\omega]^{2\ell} x_{\pm \gamma}(j) [\omega]^{-2\ell} = x_{\pm \gamma}(j \pm 2 \ell)$
	for all $\gamma\in \Gamma$. Consequently, since $e [\omega]^{-2\ell}$ commutes with the action of $\gt$ on $L(\Lambda)$, it must be equal to a scalar.
\end{dokaz}

\begin{kor} \label{omegaWgen}
	$L(\Lambda)=\langle [\omega]^2 \rangle U(\gt_1) v_\Lambda$.
\end{kor}

\bigskip

Following \cite{FS}, define extremal vectors
\begin{equation} \label{extremal vec}
 \cdots \xrightarrow{[\omega]^{-2}} v_{\Lambda} \xrightarrow{[\omega]^{-2}} v_{\Lambda}^{(-2)} \xrightarrow{[\omega]^{-2}} v_{\Lambda}^{(-4)}\xrightarrow{[\omega]^{-2}} v_{\Lambda}^{(-6)} \xrightarrow{[\omega]^{-2}}\ \cdots
\end{equation}
$$v_{\Lambda}^{(-2m)}=[\omega]^{-2m}  v_{\Lambda},\qquad m\in\Z$$ 
and the corresponding shifted Feigin-Stoyanovsky's type subspaces
$$W_{-2m}= U(\gt_1)v_{\Lambda}^{(-2m)}.$$
By Corollary \ref{omegaWgen}, 
$$L(\Lambda)=\bigcup_{m\in\Z} W_{-2m}.$$

Since
$$[\omega]^{-2m} x(\pi) v_{\Lambda}=x(\pi^{+2m}) v_{\Lambda}^{(-2m)},$$
basis of $W_{-2m}$ can be obtained by taking bases of $W$ and shifting degrees of monomials. We say that
$x(\pi)$ satisfies {\em initial conditions} $IC_{2m}$ for $W(\Lambda)$ if $x(\pi^{+2m})$ satisfies initial conditions, and that is if and only if it doesn't contain elements of degree greater than $2m$ and elements of degree $2m-1$ satisfy additional conditions corresponding to initial conditions for elements of degree $-1$ in $x(\pi)$.

\begin{prop} The set
	$$\mathcal{B}_{-2m}=\{x(\pi)v_{\Lambda}^{(-2m)} \,|\, x(\pi) \ \textrm{satisfies}\ DC\ \textrm{and}\ IC_{2m}\ \textsf{ for}\ W(\Lambda)\},$$
	is a basis for $W_{-2m}$.
\end{prop}

\bigskip

\begin{figure}[ht] \caption{Periodic path $x(\mu_r)=x(\nu_{\ell - r}^-)x(\nu_r)$.} \label{figPeriodicPath}
	\begin{center}
		\begin{tikzpicture} [scale=.17]
		
		\draw (-.5,-.5) -- +(16,0) -- +(0,16) -- cycle;
		\node at (-1,15) {$\scriptscriptstyle 1$};
		\node at (-1,14) {$\scriptscriptstyle 2$};
		\node at (-1,0) {$\scriptscriptstyle \ell$};
		\node at (0,-1) {$\scriptscriptstyle 1$};
		\node at (1,-1) {$\scriptscriptstyle 2$};
		\node at (15,-1) {$\scriptscriptstyle \ell$};
		\node[anchor=-135] at (7.5,7.5) {$\scriptscriptstyle (-1)$};
		
		\draw[dashed] (-.5,9) node[left] {$\scriptscriptstyle r$} -- (6,9);
		
		\draw[dotted] (0,8) -- (8,0);
		\foreach \n in {0,...,8}
		{
			\node[fill=red, circle, inner sep=1pt](a\n) at (\n,{8-\n}) {};
		}
		
			\node[anchor=45] at (7.5,-8.5) {$\scriptscriptstyle (-2)$};
			
			\draw (-.5,-.5) -- +(16,0) -- +(16,-16) -- cycle;
			\node at (16,-1) {$\scriptscriptstyle 1$};
			\node at (16,-2) {$\scriptscriptstyle 2$};
			\node at (16,-16) {$\scriptscriptstyle \ell$};
			
			\draw[dotted] (8.3,-.3) -- (15,-7);
			\foreach \n in {9,...,15}
			{
				\node[fill=blue, circle, inner sep=1pt](a\n) at (\n,{8-\n}) {};
			}
		\end{tikzpicture}
	\end{center}
\end{figure}
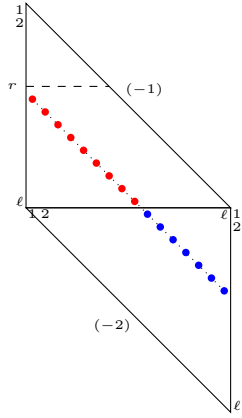

We now describe the sequence of embeddings \eqref{extremal vec} of shifted Feigin-Stoyanovsky's type subspaces $W_{-2m}$. For a level $1$ standard module $L(\Lambda_r)$ define 
$$ x(\nu_r) = x_{\ell -r,\ell} (-1) \cdots x_{2,r+2} (-1) x_{1,r+1} (-1). $$
Define also $x(\mu_r)=x(\nu_{\ell - r}^-)x(\nu_r)$ (see figure \ref{figPeriodicPath}).

\bigskip

In order to describe the action of $ x(\nu_r)$ and $x(\mu_r)$ on the highest weight vector $v_{\Lambda_r}$, we recall a few technical facts from \cite{BPT}. 
Let  $V_r=U(\g_0)v_{\Lambda_r}\subset L(\Lambda_r)$ be a $\g_0$-submodule ``at the top'' of $L(\Lambda_r)$. Then $V_r=\bigwedge^r \C^\ell$, for $r=0,\dots,\ell$, 
where $\C^\ell$ is the vector representation for $\g_0$  (cf. \cite{H}).
If $e_1,\dots,e_\ell$ is a basis for $\C^\ell$, then a basis for $V_r$ consists of vectors $v_{p_1 \dots p_r}=e_{p_1}\wedge\cdots\wedge e_{p_r}$ for $1\leq p_1<\dots<p_i\leq\ell$.
Moreover, 
\begin{equation} \label{hwvect_eq}
v_{\Lambda_r}=v_{1 2 \dots r}. 
\end{equation}
For $I=\{p_1,\dots,p_r\}$, $1\leq p_1<\dots<p_r\leq\ell$,  denote by $v_I=v_{p_1 \dots p_r}$.

The following result is from \cite{BPT} (cf. Lemmas 8, 9 and 12, Proposition 14 and Remark 3 in \cite{BPT})

\begin{prop}[{[BPT]}]\label{rel_prop}
	Let $x(\pi)=x_{i_m j_m}(-1)\cdots x_{i_2 j_2}(-1) x_{i_1 j_1}(-1)$ be such
	that $j_1\leq\dots \leq j_m$. 
	\begin{enumerate}[label=(\roman*)]
		\item \label{rel_propi} If $j_{t-1}=j_t$ and $I\subset \{1,\dots,\ell\}$ such that $\{1,\dots,j_{t-1}\}\setminus \{j_1,\dots,j_{t-1}\}\subset I$, then
		$x(\pi) v_I=0$.
		
		\item \label{rel_propii} If there is an index occuring more than twice in the sequence 
		$i_1,\dots,i_m$,\, $j_1,\dots,j_m$,\, $s_1,\dots,s_{\ell-m}$, then
		$ 
		x_{i_m j_m}(-1) \cdots x_{i_1 j_1}(-1) v_{s_1 \dots s_{\ell-m}} = 0
		$.
		
		\item \label{rel_propiii} Set $I=\{1,\dots,\ell\}\setminus \{j_1,\dots,j_m\}$, $I'=\{i_1,\dots,i_m\}$. Then
		$x(\pi) v_I =C  [\omega] v_{I'}$, $C\neq 0$.
	\end{enumerate}
\end{prop}

\bigskip

\begin{prop} \label{maxarg}
\begin{enumerate}[label=(\roman*)]
		\item \label{maxargi} $x(\nu_r)$ is the maximal monomial consisting of factors of degree greater than $-1$ that acts  nontrivially on $v_{\Lambda_r}$ and
		$$x(\nu_r)v_{\Lambda_r}=C[\omega] v_{\Lambda_{\ell - r}},\qquad C\neq 0.$$
		\item \label{maxargii}  $x(\mu_r)$ is the maximal monomial with factors of degree greater than $-2$ that acts  nontrivially on $v_{\Lambda_r}$ and
		$$x(\mu_r)v_{\Lambda_r}=C[\omega]^2 v_{\Lambda_r},\qquad C\neq 0.$$
		\item \label{maxargiii}  A monomial $x(\pi)$ satisfies difference and initial conditions for $W(\Lambda_r)$ if and only if $x(\pi)x(\mu_{r}^{+2})$ satisfies DC.
	\end{enumerate}
\end{prop}
\begin{dokaz}  For \ref{maxargi}, first note that proposition \ref{rel_prop}, \ref{rel_propiii}, and \eqref{hwvect}   immediately give
	$$x(\nu_r)v_{\Lambda_r}=C[\omega] v_{\Lambda_{\ell - r}},\qquad C\neq 0.$$
	Assume $x(\pi)=x_{i_{n} j_{n}}(-1)\cdots  x_{i_1 j_1}(-1)$ is greater than $x(\nu_r)$. Then there exists $1\leq t\leq n$ such that $i_s=s, j_s=r+s$, for $s<t$, 
	and either $t=\ell-r+1$, or $i_t=t-1$, or $i_t=t,r<j_t<r+t$. 
In the first case $$x(\pi)v_{\Lambda_r}=x_{i_{n} j_{n}}(-1)\cdots  x_{i_t j_t}(-1) C [\omega]  v_{\Lambda_{\ell - r}}= 
C [\omega] x_{i_{n} j_{n}}(0)\cdots  x_{i_t j_t}(0) v_{\Lambda_{\ell - r}} =0.$$
In the second case $x(\pi)v_{\Lambda_r}=0$ by Proposition \ref{rel_prop}, \ref{rel_propii}.
In the third case $x(\pi)v_{\Lambda_r}=0$ by Proposition \ref{rel_prop}, \ref{rel_propi}.
Claim \ref{maxargii} \quad is a direct consequence of \ref{maxargi}.
For claim \ref{maxargiii}, note that $x(\mu_{r}^{+2})=x_{r1}(0)\cdots$; the claim follows from a remark given below the definition of initial conditions.
\end{dokaz}

For a higher level module $L(\Lambda)$, $\Lambda=\Lambda_{r_1} + \dots +\Lambda_{r_k}$, with highest weight $v_{\Lambda_{r_1}} \otimes \cdots \otimes v_{\Lambda_{r_k}}$ set
	$$x(\mu_\Lambda)=x(\mu_{r_1}) \cdots x(\mu_{r_k}). $$ 
	Then Proposition \ref{maxarg} and \eqref{OrdMult_rel} give
	\begin{eqnarray*}
		x(\mu_\Lambda)v_\Lambda & = & x(\mu_{r_1}) v_{\Lambda_{r_1}} \otimes \cdots \otimes  x(\mu_{r_k})v_{\Lambda_{r_k}} \\
		& = & C_1 [\omega]^2 v_{\Lambda_{r_1}}  \otimes \cdots \otimes  C_\ell [\omega]^2 v_{\Lambda_{r_k}} \\
		& = & C [\omega]^2 v_\Lambda, 
	\end{eqnarray*}
where $ C, C_1,\dots,C_\ell\neq 0$. 

\bigskip

Note that, up to a scalar,
$$x(\pi) v_{\Lambda}	= x(\pi) [\omega]^{-2}  [\omega]^2 v_{\Lambda} = x(\pi) x(\mu_{\Lambda}^{+2}) v_{\Lambda}^{(-2)}. $$
Hence the basis of $W$ can be embedded into the basis of $W_{-2}$,  
$$\mathcal{B}_0 \subset \mathcal{B}_{-2}.$$

\begin{figure}[ht] \caption{Inclusion $\mathcal{B}_0 \subset \mathcal{B}_{-2}$: the two monomials acting on $v_{\Lambda_r}$ and $v_{\Lambda_r}^{(-2)}$}, resp., are equal.  \label{figInclusion1}
	\begin{center}
		\begin{tikzpicture} [scale=.12]
			\draw (-.5,-.5) -- +(16,0) -- +(0,16) -- +(0,0) -- +(16,-16) -- +(16,0) -- cycle;
			\draw (15.5,-16.5) -- +(16,0) -- +(0,16) -- +(0,0) -- +(1,-1);
			\draw (31.5,-17.5) -- ++(0,1) -- +(1,-1);
			\draw[dashed] (-.5,9) -- (6,9);
			
			\node[fill=orange, circle, inner sep=1pt](a) at (2,8) {};
			\node[fill=orange, circle, inner sep=1pt](b) at (5,6) {};
			\node[fill=orange, circle, inner sep=1pt](c) at (8,3) {};
			\node[fill=yellow, circle, inner sep=1pt](d) at (10,-3) {};
			\node[fill=green, circle, inner sep=1pt](e) at (16,-5) {};
			\node[fill=green, circle, inner sep=1pt](f) at (18,-9) {};
			\node[fill=green, circle, inner sep=1pt](g) at (25,-10) {};
			\draw[dotted] (a) -- (b) -- (c) -- (d) -- (e) -- (f) -- (g);
			
			\begin{scope}[xshift=-16cm,yshift=16cm]		
			\draw (-.5,-.5) -- +(16,0) -- +(0,16) -- cycle;
			
			\draw[dashed] (-.5,9)node[left] {$\scriptscriptstyle r$} -- (6,9);
			
			\draw[dotted] (0,8) -- (8,0);
			\foreach \n in {0,...,8}
			{
				\node[fill=red, circle, inner sep=1pt](a\n) at (\n,{8-\n}) {};
			}		
			
			\draw (-.5,-.5) -- +(16,0) -- +(16,-16) -- cycle;
			
			\draw[dotted] (8.3,-.3) -- (15,-7);
			\foreach \n in {9,...,15}
			{
				\node[fill=blue, circle, inner sep=1pt](a\n) at (\n,{8-\n}) {};
			}
			
			\end{scope}
			\draw[dotted] (a) -- (a15);
			\node[anchor=-135] at (-8.5,23.5) {$\scriptscriptstyle (1)$};
			\node[anchor=45] at (-8.5,7.5) {$\scriptscriptstyle (0)$};
			\node[anchor=-135] at (7.5,7.5) {$\scriptscriptstyle (-1)$};
			\node[anchor=45] at (7.5,-8.5) {$\scriptscriptstyle (-2)$};
			\node[anchor=-135] at (23.5,-8.5) {$\scriptscriptstyle (-3)$};
		
		\begin{scope}[xshift=-40cm]		
		\draw (-.5,-.5) -- +(16,0) -- +(0,16) -- +(0,0) -- +(16,-16) -- +(16,0) -- cycle;
		\draw (15.5,-16.5) -- +(16,0) -- +(0,16) -- +(0,0) -- +(1,-1);
		\draw (31.5,-17.5) -- ++(0,1) -- +(1,-1);
		\draw[dashed] (-.5,9)node[left] {$\scriptscriptstyle r$} -- (6,9) ;
		
		\node[fill=orange, circle, inner sep=1pt](xa) at (2,8) {};
		\node[fill=orange, circle, inner sep=1pt](xb) at (5,6) {};
		\node[fill=orange, circle, inner sep=1pt](xc) at (8,3) {};
		\node[fill=yellow, circle, inner sep=1pt](xd) at (10,-3) {};
		\node[fill=green, circle, inner sep=1pt](xe) at (16,-5) {};
		\node[fill=green, circle, inner sep=1pt](xf) at (18,-9) {};
		\node[fill=green, circle, inner sep=1pt](xg) at (25,-10) {};
		\draw[dotted] (xa) -- (xb) -- (xc) -- (xd) -- (xe) -- (xf) -- (xg);
		
		\node[anchor=-135] at (7.5,7.5) {$\scriptscriptstyle (-1)$};
		\node[anchor=45] at (7.5,-8.5) {$\scriptscriptstyle (-2)$};
		\node[anchor=-135] at (23.5,-8.5) {$\scriptscriptstyle (-3)$};
		\end{scope}
		\end{tikzpicture}
	\end{center}
\end{figure}
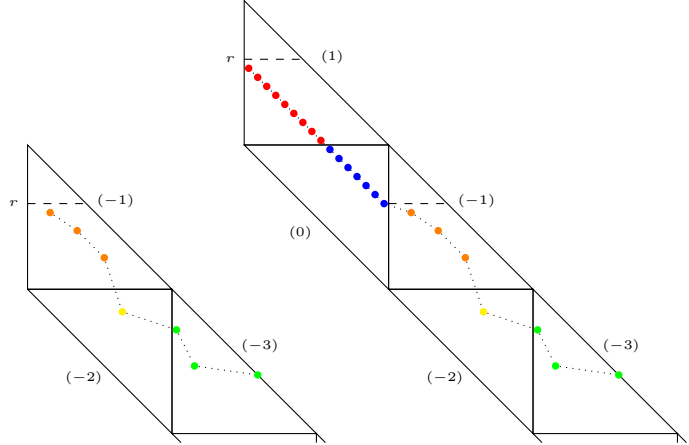

We proceed inductively and obtain a chain of inclusions
\begin{equation}\label{B_mbeddings}
\mathcal{B}_0 \subset \mathcal{B}_{-2}\subset \mathcal{B}_{-4}\subset \cdots \subset \mathcal{B}_{-2m}\subset \cdots.
\end{equation}

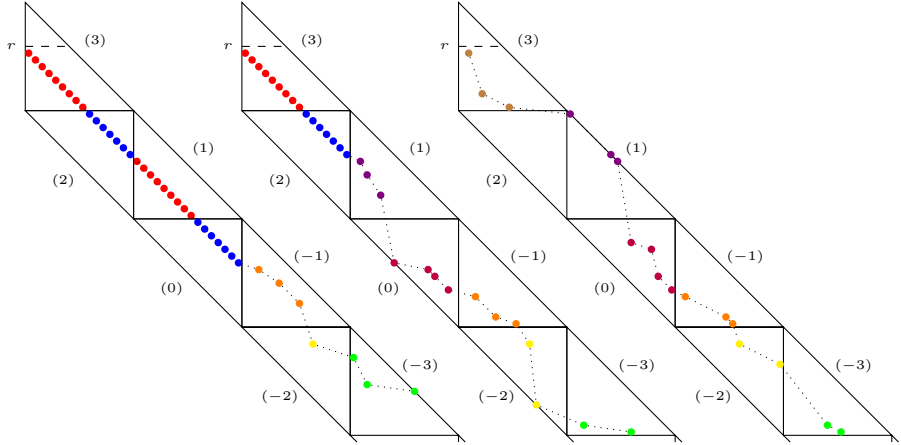
\begin{figure}[h] \caption{Inclusion of $\mathcal{B}_0 \subset \mathcal{B}_{-2}\subset \mathcal{B}_{-4}$.} \label{figInclusion2}
	\begin{center}
		\begin{tikzpicture} [scale=.09]
		\draw (-.5,-.5) -- +(16,0) -- +(0,16) -- +(0,0) -- +(16,-16) -- +(16,0) -- cycle;
		\draw (15.5,-16.5) -- +(16,0) -- +(0,16) -- +(0,0) -- +(1,-1);
		\draw (31.5,-17.5) -- ++(0,1) -- +(1,-1);
		
		\node[anchor=-135] at (-24.5,39.5) {$\scriptscriptstyle (3)$};
		\node[anchor=45] at (-24.5,23.5) {$\scriptscriptstyle (2)$};
		\node[anchor=-135] at (-8.5,23.5) {$\scriptscriptstyle (1)$};
		\node[anchor=45] at (-8.5,7.5) {$\scriptscriptstyle (0)$};
		\node[anchor=-135] at (7.5,7.5) {$\scriptscriptstyle (-1)$};
		\node[anchor=45] at (7.5,-8.5) {$\scriptscriptstyle (-2)$};
		\node[anchor=-135] at (23.5,-8.5) {$\scriptscriptstyle (-3)$};
		
		\node[fill=orange, circle, inner sep=1pt](a) at (2,8) {};
		\node[fill=orange, circle, inner sep=1pt](b) at (5,6) {};
		\node[fill=orange, circle, inner sep=1pt](c) at (8,3) {};
		\node[fill=yellow, circle, inner sep=1pt](d) at (10,-3) {};
		\node[fill=green, circle, inner sep=1pt](e) at (16,-5) {};
		\node[fill=green, circle, inner sep=1pt](f) at (18,-9) {};
		\node[fill=green, circle, inner sep=1pt](g) at (25,-10) {};
		\draw[dotted] (a) -- (b) -- (c) -- (d) -- (e) -- (f) -- (g);
		
		\begin{scope}[xshift=-16cm,yshift=16cm]		
		\draw (-.5,-.5) -- +(16,0) -- +(0,16) -- cycle;
		
		\draw[dotted] (0,8) -- (8,0);
		\foreach \n in {0,...,8} { \node[fill=red, circle, inner sep=1pt](a\n) at (\n,{8-\n}) {}; }
		
		\draw (-.5,-.5) -- +(16,0) -- +(16,-16) -- cycle;
		
		\draw[dotted] (8.3,-.3) -- (15,-7);
		\foreach \n in {9,...,15} {\node[fill=blue, circle, inner sep=1pt](a\n) at (\n,{8-\n}) {};}		
		\end{scope}
		\draw[dotted] (a) -- (a15);
		
		\begin{scope}[xshift=-32cm,yshift=32cm]		
		\draw (-.5,-.5) -- +(16,0) -- +(0,16) -- cycle;
		\draw[dashed] (-.5,9)  node[left] {$\scriptscriptstyle r$}-- (6,9);
		
		\draw[dotted] (0,8) -- (8,0);
		\foreach \n in {0,...,8} { \node[fill=red, circle, inner sep=1pt](b\n) at (\n,{8-\n}) {}; }
		
		\draw (-.5,-.5) -- +(16,0) -- +(16,-16) -- cycle;
		
		\draw[dotted] (8.3,-.3) -- (15,-7);
		\foreach \n in {9,...,15} {\node[fill=blue, circle, inner sep=1pt](b\n) at (\n,{8-\n}) {};}		
		\end{scope}
		\draw[dotted] (b15) -- (a1);

		\begin{scope} [xshift=32cm]
		\draw (-.5,-.5) -- +(16,0) -- +(0,16) -- +(0,0) -- +(16,-16) -- +(16,0) -- cycle;
		\draw (15.5,-16.5) -- +(16,0) -- +(0,16) -- +(0,0) -- +(1,-1);
		\draw (31.5,-17.5) -- ++(0,1) -- +(1,-1);
		
		\node[anchor=-135] at (-24.5,39.5) {$\scriptscriptstyle (3)$};
		\node[anchor=45] at (-24.5,23.5) {$\scriptscriptstyle (2)$};
		\node[anchor=-135] at (-8.5,23.5) {$\scriptscriptstyle (1)$};
		\node[anchor=45] at (-8.5,7.5) {$\scriptscriptstyle (0)$};
		\node[anchor=-135] at (7.5,7.5) {$\scriptscriptstyle (-1)$};
		\node[anchor=45] at (7.5,-8.5) {$\scriptscriptstyle (-2)$};
		\node[anchor=-135] at (23.5,-8.5) {$\scriptscriptstyle (-3)$};
		
		\node[fill=orange, circle, inner sep=1pt](a) at (2,4) {};
		\node[fill=orange, circle, inner sep=1pt](b) at (5,1) {};
		\node[fill=orange, circle, inner sep=1pt](c) at (8,0) {};
		\node[fill=yellow, circle, inner sep=1pt](d) at (10,-3) {};
		\node[fill=yellow, circle, inner sep=1pt](e) at (11,-12) {};
		\node[fill=green, circle, inner sep=1pt](f) at (18,-15) {};
		\node[fill=green, circle, inner sep=1pt](g) at (25,-16) {};
		\draw[dotted] (a) -- (b) -- (c) -- (d) -- (e) -- (f) -- (g);
		
		\begin{scope}[xshift=-16cm,yshift=16cm]		
		\draw (-.5,-.5) -- +(16,0) -- +(0,16) -- cycle;
		\draw (-.5,-.5) -- +(16,0) -- +(16,-16) -- cycle;
		
		\node[fill=violet, circle, inner sep=1pt](aa1) at (1,8) {};
		\node[fill=violet, circle, inner sep=1pt](bb1) at (2,6) {};
		\node[fill=violet, circle, inner sep=1pt](c1) at (4,3) {};
		\node[fill=purple, circle, inner sep=1pt](d1) at (6,-7) {};
		\node[fill=purple, circle, inner sep=1pt](e1) at (11,-8) {};
		\node[fill=purple, circle, inner sep=1pt](f1) at (12,-9) {};
		\node[fill=purple, circle, inner sep=1pt](g1) at (14,-11) {};
		\draw[dotted] (aa1) -- (bb1) -- (c1) -- (d1) -- (e1) -- (f1) -- (g1);
		\end{scope}
		\draw[dotted] (a) -- (g1);
		
		\begin{scope}[xshift=-32cm,yshift=32cm]		
		\draw (-.5,-.5) -- +(16,0) -- +(0,16) -- cycle;
		\draw[dashed] (-.5,9)node[left] {$\scriptscriptstyle r$} -- (6,9) ;
		
		\draw[dotted] (0,8) -- (8,0);
		\foreach \n in {0,...,8} { \node[fill=red, circle, inner sep=1pt](b\n) at (\n,{8-\n}) {}; }
		
		\draw (-.5,-.5) -- +(16,0) -- +(16,-16) -- cycle;
		
		\draw[dotted] (8.3,-.3) -- (15,-7);
		\foreach \n in {9,...,15} {\node[fill=blue, circle, inner sep=1pt](b\n) at (\n,{8-\n}) {};}		
		\end{scope}
		\draw[dotted] (b15) -- (aa1);
		
		\end{scope}
		
		\begin{scope} [xshift=64cm]
		\draw (-.5,-.5) -- +(16,0) -- +(0,16) -- +(0,0) -- +(16,-16) -- +(16,0) -- cycle;
		\draw (15.5,-16.5) -- +(16,0) -- +(0,16) -- +(0,0) -- +(1,-1);
		\draw (31.5,-17.5) -- ++(0,1) -- +(1,-1);
		
		\node[anchor=-135] at (-24.5,39.5) {$\scriptscriptstyle (3)$};
		\node[anchor=45] at (-24.5,23.5) {$\scriptscriptstyle (2)$};
		\node[anchor=-135] at (-8.5,23.5) {$\scriptscriptstyle (1)$};
		\node[anchor=45] at (-8.5,7.5) {$\scriptscriptstyle (0)$};
		\node[anchor=-135] at (7.5,7.5) {$\scriptscriptstyle (-1)$};
		\node[anchor=45] at (7.5,-8.5) {$\scriptscriptstyle (-2)$};
		\node[anchor=-135] at (23.5,-8.5) {$\scriptscriptstyle (-3)$};
		
		\node[fill=orange, circle, inner sep=1pt](a) at (1,4) {};
		\node[fill=orange, circle, inner sep=1pt](b) at (7,1) {};
		\node[fill=orange, circle, inner sep=1pt](c) at (8,0) {};
		\node[fill=yellow, circle, inner sep=1pt](d) at (9,-3) {};
		\node[fill=yellow, circle, inner sep=1pt](e) at (15,-6) {};
		\node[fill=green, circle, inner sep=1pt](f) at (22,-15) {};
		\node[fill=green, circle, inner sep=1pt](g) at (24,-16) {};
		\draw[dotted] (a) -- (b) -- (c) -- (d) -- (e) -- (f) -- (g);
		
		\begin{scope}[xshift=-16cm,yshift=16cm]		
		\draw (-.5,-.5) -- +(16,0) -- +(0,16) -- cycle;
		\draw (-.5,-.5) -- +(16,0) -- +(16,-16) -- cycle;
		
		\node[fill=violet, circle, inner sep=1pt](aa1) at (0,15) {};
		\node[fill=violet, circle, inner sep=1pt](bb1) at (6,9) {};
		\node[fill=violet, circle, inner sep=1pt](c1) at (7,8) {};
		\node[fill=purple, circle, inner sep=1pt](d1) at (9,-4) {};
		\node[fill=purple, circle, inner sep=1pt](e1) at (12,-5) {};
		\node[fill=purple, circle, inner sep=1pt](f1) at (13,-9) {};
		\node[fill=purple, circle, inner sep=1pt](g1) at (15,-11) {};
		\draw[dotted] (aa1) -- (bb1) -- (c1) -- (d1) -- (e1) -- (f1) -- (g1);
		\end{scope}
		\draw[dotted] (a) -- (g1);
		
		\begin{scope}[xshift=-32cm,yshift=32cm]		
		\draw (-.5,-.5) -- +(16,0) -- +(0,16) -- (-.5,-.5) -- +(16,-16);
		\draw[dashed] (-.5,9)node[left] {$\scriptscriptstyle r$} -- (6,9) ;
		\node[fill=brown, circle, inner sep=1pt](aaa1) at (1,8) {};
		\node[fill=brown, circle, inner sep=1pt](bbb1) at (3,2) {};
		\node[fill=brown, circle, inner sep=1pt](ccc1) at (7,0) {};
		
		\end{scope}
		\draw[dotted] (aaa1) -- (bbb1) -- (ccc1) -- (aa1);
		
		\end{scope}
		
		\end{tikzpicture}
	\end{center}
\end{figure}

Now take ``the limit $m\to\infty$'', i.e. take inductive limit of the sequence \eqref{B_mbeddings}, to obtain a basis of the whole $L(\Lambda)$. 
Formally, set 
$$v_{\Lambda}^{(-2m)} = x(\mu_{\Lambda}^{+2m+2}) x(\mu_{\Lambda}^{+2m+4})\cdots v_{\Lambda}^{(-\infty)}.$$ 
We say that a semi-infinite monomial $x(\pi)$ {\em has a periodic tail}, or that it {\em stabilizes}, if from some point on it consists of successive shifts of $x(\mu_\Lambda)$,
$$x(\pi)=x(\pi') x(\mu_{\Lambda}^{+2m+2}) x(\mu_{\Lambda}^{+2m+4})\cdots,$$
for some $m\in\Z$.

\begin{tm} The set
		\begin{equation*}
			\overline{\mathcal{B}}  =  
			\{ x(\pi) v_{\Lambda}^{(-\infty)} \,|\, 
			x(\pi) \ \textrm{stabilizes and satisfies}\ DC\}
		\end{equation*}
		is a basis for $L(\Lambda)$.
\end{tm}

\begin{figure}[h] \caption{Basis of standard modules consisting of semi-infinite monomials  that stabilize and satisfy DC
		} \label{figBasisL}
	\begin{center}
		\begin{tikzpicture} [scale=.1]
		\draw (-.5,-.5) -- +(16,0) -- +(0,16) -- +(0,0) -- +(16,-16) -- +(16,0) -- cycle;
		\draw (15.5,-16.5) -- +(16,0) -- +(0,16) -- +(0,0) -- +(1,-1);
		\draw (31.5,-17.5) -- ++(0,1) -- +(1,-1);
		
		\node[anchor=-135] at (-40.5,55.5) {$\scriptscriptstyle (5)$};
		\node[anchor=45] at (-40.5,39.5) {$\scriptscriptstyle (4)$};
		\node[anchor=-135] at (-24.5,39.5) {$\scriptscriptstyle (3)$};
		\node[anchor=45] at (-24.5,23.5) {$\scriptscriptstyle (2)$};
		\node[anchor=-135] at (-8.5,23.5) {$\scriptscriptstyle (1)$};
		\node[anchor=45] at (-8.5,7.5) {$\scriptscriptstyle (0)$};
		\node[anchor=-135] at (7.5,7.5) {$\scriptscriptstyle (-1)$};
		\node[anchor=45] at (7.5,-8.5) {$\scriptscriptstyle (-2)$};
		\node[anchor=-135] at (23.5,-8.5) {$\scriptscriptstyle (-3)$};
		
		\node[fill=orange, circle, inner sep=1pt](a) at (2,8) {};
		\node[fill=orange, circle, inner sep=1pt](b) at (5,6) {};
		\node[fill=orange, circle, inner sep=1pt](c) at (8,3) {};
		\node[fill=yellow, circle, inner sep=1pt](d) at (10,-3) {};
		\node[fill=green, circle, inner sep=1pt](e) at (16,-5) {};
		\node[fill=green, circle, inner sep=1pt](f) at (18,-9) {};
		\node[fill=green, circle, inner sep=1pt](g) at (25,-10) {};
		\draw[dotted] (a) -- (b) -- (c) -- (d) -- (e) -- (f) -- (g);
		
		\begin{scope}[xshift=-16cm,yshift=16cm]		
		\draw (-.5,-.5) -- +(16,0) -- +(0,16) -- cycle;
		
		\draw[dotted] (0,8) -- (8,0);
		\foreach \n in {0,...,8} { \node[fill=red, circle, inner sep=1pt](a\n) at (\n,{8-\n}) {}; }
		
		\draw (-.5,-.5) -- +(16,0) -- +(16,-16) -- cycle;
		
		\draw[dotted] (8.3,-.3) -- (15,-7);
		\foreach \n in {9,...,15} {\node[fill=blue, circle, inner sep=1pt](a\n) at (\n,{8-\n}) {};}		
		\end{scope}
		\draw[dotted] (a) -- (a15);
		
		\begin{scope}[xshift=-32cm,yshift=32cm]		
		\draw (-.5,-.5) -- +(16,0) -- +(0,16) -- cycle;
		
		\draw[dotted] (0,8) -- (8,0);
		\foreach \n in {0,...,8} { \node[fill=red, circle, inner sep=1pt](b\n) at (\n,{8-\n}) {}; }
		
		\draw (-.5,-.5) -- +(16,0) -- +(16,-16) -- cycle;
		
		\draw[dotted] (8.3,-.3) -- (15,-7);
		\foreach \n in {9,...,15} {\node[fill=blue, circle, inner sep=1pt](b\n) at (\n,{8-\n}) {};}		
		\end{scope}
		\draw[dotted] (b15) -- (a1);
		
		\begin{scope}[xshift=-48cm,yshift=48cm]		
		\draw (-.5,-.5) -- +(16,0) -- +(0,16) -- cycle;
		
		\draw[dotted] (0,8) -- (8,0);
		\foreach \n in {0,...,8} { \node[fill=red, circle, inner sep=1pt](h\n) at (\n,{8-\n}) {}; }
		
		\draw (-.5,-.5) -- +(16,0) -- +(16,-16) -- cycle;
		
		\draw[dotted] (8.3,-.3) -- (15,-7);
		\foreach \n in {9,...,15} {\node[fill=blue, circle, inner sep=1pt](h\n) at (\n,{8-\n}) {};}		
		
		\node[fill=blue, circle, inner sep=1pt] at (-1,9) {};
		\end{scope}
		\draw[dotted] (h15) -- (b1);
		
		\draw (-48.5,47.5) -- +(-1,1); 
		\draw (-49.5,63.5) -- ++(1,0) -- +(-1,1); 
		
		\begin{scope} [xshift=35cm]
		\draw (-.5,-.5) -- +(16,0) -- +(0,16) -- +(0,0) -- +(16,-16) -- +(16,0) -- cycle;
		\draw (15.5,-16.5) -- +(16,0) -- +(0,16) -- +(0,0) -- +(1,-1);
		\draw (31.5,-17.5) -- ++(0,1) -- +(1,-1);
		
		\node[anchor=-135] at (-40.5,55.5) {$\scriptscriptstyle (5)$};
		\node[anchor=45] at (-40.5,39.5) {$\scriptscriptstyle (4)$};
		\node[anchor=-135] at (-24.5,39.5) {$\scriptscriptstyle (3)$};
		\node[anchor=45] at (-24.5,23.5) {$\scriptscriptstyle (2)$};
		\node[anchor=-135] at (-8.5,23.5) {$\scriptscriptstyle (1)$};
		\node[anchor=45] at (-8.5,7.5) {$\scriptscriptstyle (0)$};
		\node[anchor=-135] at (7.5,7.5) {$\scriptscriptstyle (-1)$};
		\node[anchor=45] at (7.5,-8.5) {$\scriptscriptstyle (-2)$};
		\node[anchor=-135] at (23.5,-8.5) {$\scriptscriptstyle (-3)$};
		
		\node[fill=orange, circle, inner sep=1pt](a) at (2,4) {};
		\node[fill=orange, circle, inner sep=1pt](b) at (5,1) {};
		\node[fill=orange, circle, inner sep=1pt](c) at (8,0) {};
		\node[fill=yellow, circle, inner sep=1pt](d) at (10,-3) {};
		\node[fill=yellow, circle, inner sep=1pt](e) at (11,-12) {};
		\node[fill=green, circle, inner sep=1pt](f) at (18,-15) {};
		\node[fill=green, circle, inner sep=1pt](g) at (25,-16) {};
		\draw[dotted] (a) -- (b) -- (c) -- (d) -- (e) -- (f) -- (g);
		
		\begin{scope}[xshift=-16cm,yshift=16cm]		
		\draw (-.5,-.5) -- +(16,0) -- +(0,16) -- cycle;
		\draw (-.5,-.5) -- +(16,0) -- +(16,-16) -- cycle;
		
		\node[fill=violet, circle, inner sep=1pt](aa1) at (1,8) {};
		\node[fill=violet, circle, inner sep=1pt](bb1) at (2,6) {};
		\node[fill=violet, circle, inner sep=1pt](c1) at (4,3) {};
		\node[fill=purple, circle, inner sep=1pt](d1) at (6,-7) {};
		\node[fill=purple, circle, inner sep=1pt](e1) at (11,-8) {};
		\node[fill=purple, circle, inner sep=1pt](f1) at (12,-9) {};
		\node[fill=purple, circle, inner sep=1pt](g1) at (14,-11) {};
		\draw[dotted] (aa1) -- (bb1) -- (c1) -- (d1) -- (e1) -- (f1) -- (g1);
		\end{scope}
		\draw[dotted] (a) -- (g1);
		
		\begin{scope}[xshift=-32cm,yshift=32cm]		
		\draw (-.5,-.5) -- +(16,0) -- +(0,16) -- cycle;
		
		\draw[dotted] (0,8) -- (8,0);
		\foreach \n in {0,...,8} { \node[fill=red, circle, inner sep=1pt](b\n) at (\n,{8-\n}) {}; }
		
		\draw (-.5,-.5) -- +(16,0) -- +(16,-16) -- cycle;
		
		\draw[dotted] (8.3,-.3) -- (15,-7);
		\foreach \n in {9,...,15} {\node[fill=blue, circle, inner sep=1pt](b\n) at (\n,{8-\n}) {};}		
		\end{scope}
		\draw[dotted] (b15) -- (aa1);
		
		\begin{scope}[xshift=-48cm,yshift=48cm]		
		\draw (-.5,-.5) -- +(16,0) -- +(0,16) -- cycle;
		
		\draw[dotted] (0,8) -- (8,0);
		\foreach \n in {0,...,8} { \node[fill=red, circle, inner sep=1pt](h\n) at (\n,{8-\n}) {}; }
		
		\draw (-.5,-.5) -- +(16,0) -- +(16,-16) -- cycle;
		
		\draw[dotted] (8.3,-.3) -- (15,-7);
		\foreach \n in {9,...,15} {\node[fill=blue, circle, inner sep=1pt](h\n) at (\n,{8-\n}) {};}		
		
		\node[fill=blue, circle, inner sep=1pt] at (-1,9) {};
		\end{scope}
		\draw[dotted] (h15) -- (b1);

		\draw (-48.5,47.5) -- +(-1,1); 
		\draw (-49.5,63.5) -- ++(1,0) -- +(-1,1); 
		\end{scope}
		
		\end{tikzpicture}
	\end{center}
\end{figure}
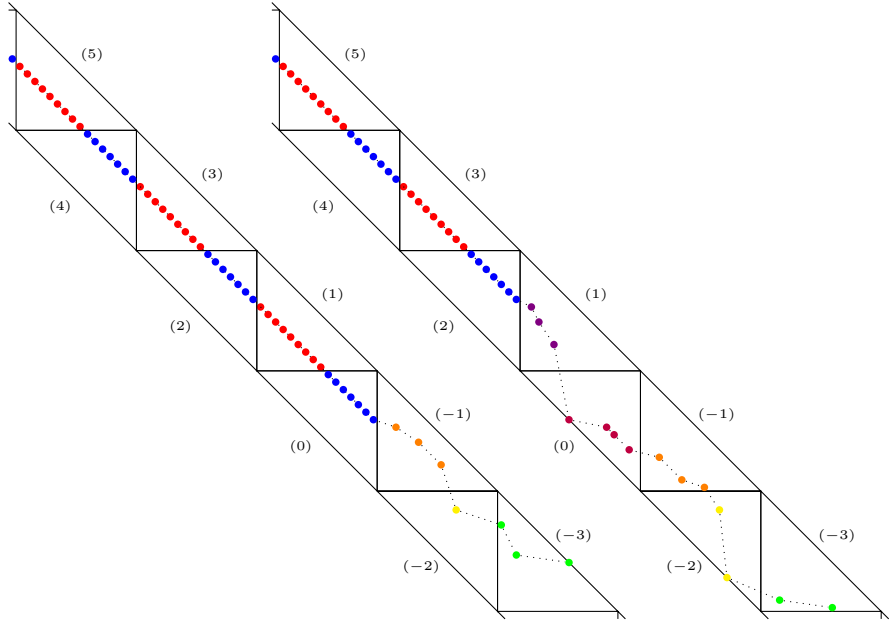

\end{document}